\documentclass[11pt, a4paper]{article}
\usepackage[cp1251]{inputenc}
\usepackage[russian]{babel}
\usepackage{amsmath} \usepackage{euscript}
\usepackage{longtable}
\usepackage{mathrsfs}
\usepackage{amssymb}
\begin{document}

\newcounter{bnomer} \newcounter{snomer}
\newcounter{bsnomer}
\setcounter{bnomer}{0}
\renewcommand{\thesnomer}{\thebnomer.\arabic{snomer}}
\renewcommand{\thebsnomer}{\thebnomer.\arabic{bsnomer}}
\renewcommand{\refname}{\begin{center}\large{\textbf{References}}\end{center}}

\newcommand{\sect}[1]{%
\setcounter{snomer}{0}\setcounter{bsnomer}{0}
\refstepcounter{bnomer}
\par\bigskip\begin{center}\large{\textbf{\arabic{bnomer}. {#1}}}\end{center}}
\newcommand{\sst}{%
\refstepcounter{bsnomer}
\par\bigskip\textbf{\arabic{bnomer}.\arabic{bsnomer}. }}
\newcommand{\defi}[1]{%
\refstepcounter{snomer}
\par\medskip\textbf{Definition \arabic{bnomer}.\arabic{snomer}. }{#1}\par\medskip}
\newcommand{\theo}[2]{%
\refstepcounter{snomer}
\par\textbf{Теорема \arabic{bnomer}.\arabic{snomer}. }{#2} {\emph{#1}}\hspace{\fill}$\square$\par}
\newcommand{\mtheop}[2]{%
\refstepcounter{snomer}
\par\textbf{Theorem \arabic{bnomer}.\arabic{snomer}. }{\emph{#1}}
\par\textsc{Proof}. {#2}\hspace{\fill}$\square$\par}
\newcommand{\mcorop}[2]{%
\refstepcounter{snomer}
\par\textbf{Corollary \arabic{bnomer}.\arabic{snomer}. }{\emph{#1}}
\par\textsc{Proof}. {#2}\hspace{\fill}$\square$\par}
\newcommand{\mtheo}[1]{%
\refstepcounter{snomer}
\par\medskip\textbf{Theorem \arabic{bnomer}.\arabic{snomer}. }{\emph{#1}}\par\medskip}
\newcommand{\mlemm}[1]{%
\refstepcounter{snomer}
\par\medskip\textbf{Lemma \arabic{bnomer}.\arabic{snomer}. }{\emph{#1}}\par\medskip}
\newcommand{\mprop}[1]{%
\refstepcounter{snomer}
\par\medskip\textbf{Proposition \arabic{bnomer}.\arabic{snomer}. }{\emph{#1}}\par\medskip}
\newcommand{\theobp}[2]{%
\refstepcounter{snomer}
\par\textbf{Теорема \arabic{bnomer}.\arabic{snomer}. }{#2} {\emph{#1}}\par}
\newcommand{\theop}[2]{%
\refstepcounter{snomer}
\par\textbf{Theorem \arabic{bnomer}.\arabic{snomer}. }{\emph{#1}}
\par\textsc{Proof}. {#2}\hspace{\fill}$\square$\par}
\newcommand{\theosp}[2]{%
\refstepcounter{snomer}
\par\textbf{Теорема \arabic{bnomer}.\arabic{snomer}. }{\emph{#1}}
\par\textbf{Схема доказательства}. {#2}\hspace{\fill}$\square$\par}
\newcommand{\exam}[1]{%
\refstepcounter{snomer}
\par\medskip\textbf{Example \arabic{bnomer}.\arabic{snomer}. }{#1}\par\medskip}
\newcommand{\deno}[1]{%
\refstepcounter{snomer}
\par\textbf{Definition \arabic{bnomer}.\arabic{snomer}. }{#1}\par}
\newcommand{\post}[1]{%
\refstepcounter{snomer}
\par\textbf{Предложение \arabic{bnomer}.\arabic{snomer}. }{\emph{#1}}\hspace{\fill}$\square$\par}
\newcommand{\postp}[2]{%
\refstepcounter{snomer}
\par\medskip\textbf{Proposition \arabic{bnomer}.\arabic{snomer}. }{\emph{#1}}%
\ifhmode\par\fi\textsc{Proof}. {#2}\hspace{\fill}$\square$\par\medskip}
\newcommand{\lemm}[1]{%
\refstepcounter{snomer}
\par\textbf{Lemma \arabic{bnomer}.\arabic{snomer}. }{\emph{#1}}\hspace{\fill}$\square$\par}
\newcommand{\lemmp}[2]{%
\refstepcounter{snomer}
\par\medskip\textbf{Lemma \arabic{bnomer}.\arabic{snomer}. }{\emph{#1}}
\par\textsc{Proof}. {#2}\hspace{\fill}$\square$\par\medskip}
\newcommand{\coro}[1]{%
\refstepcounter{snomer}
\par\textbf{Следствие \arabic{bnomer}.\arabic{snomer}. }{\emph{#1}}\hspace{\fill}$\square$\par}
\newcommand{\mcoro}[1]{%
\refstepcounter{snomer}
\par\textbf{Corollary \arabic{bnomer}.\arabic{snomer}. }{\emph{#1}}\par\medskip}
\newcommand{\corop}[2]{%
\refstepcounter{snomer}
\par\textbf{Следствие \arabic{bnomer}.\arabic{snomer}. }{\emph{#1}}
\par\textsc{Proof}. {#2}\hspace{\fill}$\square$\par}
\newcommand{\nota}[1]{%
\refstepcounter{snomer}
\par\medskip\textbf{Remark \arabic{bnomer}.\arabic{snomer}. }{#1}\par\medskip}
\newcommand{\propp}[2]{%
\refstepcounter{snomer}
\par\medskip\textbf{Proposition \arabic{bnomer}.\arabic{snomer}. }{\emph{#1}}
\par\textsc{Proof}. {#2}\hspace{\fill}$\square$\par\medskip}
\newcommand{\hypo}[1]{%
\refstepcounter{snomer}
\par\medskip\textbf{Conjecture \arabic{bnomer}.\arabic{snomer}. }{\emph{#1}}\par\medskip}
\newcommand{\prop}[1]{%
\refstepcounter{snomer}
\par\textbf{Proposition \arabic{bnomer}.\arabic{snomer}. }{\emph{#1}}\hspace{\fill}$\square$\par}

\newcommand{\Ind}[3]{%
\mathrm{Ind}_{#1}^{#2}{#3}}
\newcommand{\Res}[3]{%
\mathrm{Res}_{#1}^{#2}{#3}}
\newcommand{\epsi}{\epsilon}
\newcommand{\tri}{\triangleleft}
\newcommand{\Supp}[1]{%
\mathrm{Supp}(#1)}

\newcommand{\reg}{\mathrm{reg}}
\newcommand{\empr}[2]{[-{#1},{#1}]\times[-{#2},{#2}]}
\newcommand{\sreg}{\mathrm{sreg}}
\newcommand{\codim}{\mathrm{codim}\,}
\newcommand{\chara}{\mathrm{char}\,}
\newcommand{\rk}{\mathrm{rk}\,}
\newcommand{\chr}{\mathrm{ch}\,}
\newcommand{\id}{\mathrm{id}}
\newcommand{\Ad}{\mathrm{Ad}}
\newcommand{\col}{\mathrm{col}}
\newcommand{\row}{\mathrm{row}}
\newcommand{\low}{\mathrm{low}}
\newcommand{\pho}{\hphantom{\quad}\vphantom{\mid}}
\newcommand{\fho}[1]{\vphantom{\mid}\setbox0\hbox{00}\hbox to \wd0{\hss\ensuremath{#1}\hss}}
\newcommand{\wt}{\widetilde}
\newcommand{\wh}{\widehat}
\newcommand{\ad}[1]{\mathrm{ad}_{#1}}
\newcommand{\tr}{\mathrm{tr}\,}
\newcommand{\GL}{\mathrm{GL}}
\newcommand{\SL}{\mathrm{SL}}
\newcommand{\SO}{\mathrm{SO}}
\newcommand{\Sp}{\mathrm{Sp}}
\newcommand{\Mat}{\mathrm{Mat}}
\newcommand{\Stab}{\mathrm{Stab}}

\newcommand{\vfi}{\varphi}
\newcommand{\teta}{\vartheta}
\newcommand{\Bfi}{\Phi}
\newcommand{\Fp}{\mathbb{F}}
\newcommand{\Rp}{\mathbb{R}}
\newcommand{\Zp}{\mathbb{Z}}
\newcommand{\Cp}{\mathbb{C}}
\newcommand{\ut}{\mathfrak{u}}
\newcommand{\at}{\mathfrak{a}}
\newcommand{\nt}{\mathfrak{n}}
\newcommand{\mt}{\mathfrak{m}}
\newcommand{\htt}{\mathfrak{h}}
\newcommand{\spt}{\mathfrak{sp}}
\newcommand{\rt}{\mathfrak{r}}
\newcommand{\rad}{\mathfrak{rad}}
\newcommand{\bt}{\mathfrak{b}}
\newcommand{\gt}{\mathfrak{g}}
\newcommand{\vt}{\mathfrak{v}}
\newcommand{\pt}{\mathfrak{p}}
\newcommand{\Xt}{\mathfrak{X}}
\newcommand{\Po}{\mathcal{P}}
\newcommand{\Uo}{\EuScript{U}}
\newcommand{\Fo}{\EuScript{F}}
\newcommand{\Do}{\EuScript{D}}
\newcommand{\Eo}{\EuScript{E}}
\newcommand{\Iu}{\mathcal{I}}
\newcommand{\Mo}{\mathcal{M}}
\newcommand{\Nu}{\mathcal{N}}
\newcommand{\Ro}{\mathcal{R}}
\newcommand{\Co}{\mathcal{C}}
\newcommand{\Lo}{\mathcal{L}}
\newcommand{\Ou}{\mathcal{O}}
\newcommand{\Uu}{\mathcal{U}}
\newcommand{\Au}{\mathcal{A}}
\newcommand{\Vu}{\mathcal{V}}
\newcommand{\Bu}{\mathcal{B}}
\newcommand{\Sy}{\mathcal{Z}}
\newcommand{\Sb}{\mathcal{F}}
\newcommand{\Gr}{\mathcal{G}}
\newcommand{\rtc}[1]{C_{#1}^{\mathrm{red}}}

\author{Mikhail V. Ignatyev\and Aleksandr A. Shevchenko}

\date{Samara National Research University\\\texttt{mihail.ignatev@gmail.com}\hspace{1cm}\texttt{shevchenko.alexander.1618@gmail.com}}
\title{\Large{On tangent cones to Schubert varieties in type $E$}\mbox{$\vphantom{1}$}\footnotetext{The authors were supported by the Foundation for the Advancement of Theoretical Physics and Mathematics ``BASIS'', grant no. 18--1--7--2--1. The first author was also partially supported by RFBR grant no. 20--01--00091a.}} \maketitle

\begin{center}
\begin{tabular}{p{15cm}}
\small{\textsc{Abstract}. We consider tangent cones to Schubert subvarieties of the flag variety $G/B$, where $B$ is a Borel subgroup of a reductive complex algebraic group $G$ of type $E_6$, $E_7$ or $E_8$. We prove that if $w_1$ and $w_2$ form a good pair of involutions in the Weyl group $W$ of $G$ then the tangent cones $C_{w_1}$ and $C_{w_2}$ to the corresponding Schubert subvarieties of $G/B$ do not coincide as subschemes of the tangent space to $G/B$ at the neutral point.}
\end{tabular}
\end{center}

\sect{Introduction and the main result}

\sst Let $G$ be a complex reductive algebraic group, $T$ a maximal torus in~$G$, $B$ a Borel subgroup in~$G$ containing $T$, and $U$ the unipotent radical of $B$. Let $\Phi$ be the root system of $G$ with respect to~$T$, $\Phi^+$ the set of positive roots with respect to $B$, $\Delta$ the set of simple roots, and $W$ the Weyl group of $\Phi$ (see \cite{Bourbaki}, \cite{Humphreys} and \cite{Humpreys2} for basic facts about algebraic groups and root systems).

Denote by $\Fo=G/B$ the flag variety and by $X_w\subseteq\Fo$ the Schubert subvariety corresponding to an element $w$ of the Weyl group $W$. Denote by $\Ou=\Ou_{p,X_w}$ the local ring at the point $p=eB\in X_w$. Let $\mt$ be the maximal ideal of~$\Ou$. The decreasing sequence of ideals $$\Ou\supseteq\mt\supseteq\mt^2\supseteq\ldots$$ is a filtration on $\Ou$. We define $R$ to be the graded algebra $$R=\mathrm{gr}\,\Ou=\bigoplus_{i\geq0}\mt^i/\mt^{i+1}.$$ By definition, the \emph{tangent cone} $C_w$ to the Schubert variety $X_w$ at the point $p$ is the spectrum of~$R$: $C_w=\mathrm{Spec}\,R$. Obviously, $C_w$ is a subscheme of the tangent space $T_pX_w\subseteq T_p\Fo$. A hard problem in studying geometry of $X_w$ is to describe $C_w$ \cite[Chapter 7]{BilleyLakshmibai}.

In 2011, D.Yu. Eliseev and A.N. Panov computed tangent cones $C_w$ for all $w\in W$ in the case $G=\mathrm{SL}_n(\mathbb{C})$, $n\leq5$ \cite{EliseevPanov}. Using their computations, A.N. Panov formulated the following Conjecture.

\hypo{\textup{(A.N. Panov, 2011)} Let $w_1$\textup{,} $w_2$ be \label{mconj}involutions\textup{,} i.e.\textup{,} $w_1^2=w_2^2=\id$. If $w_1\neq w_2$\textup{,} then $C_{w_1}\neq C_{w_2}$ as subschemes of $T_p\Fo$.}

One can easily check that it is enough to prove the Conjecture for irreducible root systems (see Remark~\ref{nota:irred} below). In 2013, D.Yu. Eliseev and the first author proved this Conjecture in types $A_n$, $F_4$ and $G_2$ \cite{EliseevIgnatyev}. In \cite{BochkarevIgnatyevShevchenko}, M.A. Bochkarev and the authors proved the Conjecture in types $B_n$ and~$C_n$. In \cite{IgnatyevShevchenko1}, we proved that the Conjecture is true if $\Phi$ is of type $D_n$ and $w_1$, $w_2$ are so-called basic involutions. In this paper, we prove that the Conjecture is true for so-called good pairs of involutions (see Definition~\ref{defi:good_pair}) for $\Phi=E_6$, $E_7$ and $E_8$. Precisely, our main result is as follows.

\mtheo{Assume that every irreducible component of $\Phi$ is of type $E_6$, $E_7$ or $E_8$. Let $w_1$\textup,~$w_2$ be a good pair of involutions in the Weyl group of $\Phi$. Then the tangent cones $C_{w_1}$ and~$C_{w_2}$ do~not coincide as subschemes of $T_p\Fo$.\label{mtheo:non_red}}

\nota{One can also consider reduced tangent cones. Let $\Au$ be the symmetric algebra of the vector space $\mt/\mt^2$, or, equivalently, the algebra of regular functions on the tangent space $T_pX_w$. Since $R$ is generated as $\Cp$-algebra by $\mt/\mt^2$, it is a~quotient ring $R=\Au/I$. By definition, the \emph{reduced tangent cone} $C_w^{\mathrm{red}}$ to $X_w$ at the point $p$ is the common zero locus in $T_pX_w$ of the polynomials $f\in I\subseteq\Au$. Clearly, if $\rtc{w_1}\neq\rtc{w_2}$, then $C_{w_1}\neq C_{w_2}$. It was proved in \cite{BochkarevIgnatyevShevchenko} that if $\Phi$ is of type $B_n$ or $C_n$ and $w_1$ and $w_2$ are distinct involutions in $W$, then $\rtc{w_1}$ and $\rtc{w_2}$ do not coincide as subvarieties of $T_p\Fo$. In \cite{IgnatyevShevchenko1}, the similar result was obtained for basic involutions in type $D_n$. For type $E$, this question still remains open even for good pairs of involutions.}

The paper is organized as follows. In Section~\ref{sect:Kostant_Kumar}, we introduce the main technical tool used in the proof of Theorem~\ref{mtheo:non_red}. Namely, to each element $w\in W$ one can assign a polynomial $d_w$ in the algebra of regular functions on the Lie algebra of the maximal torus $T$. These polynomials are called Kostant--Kumar polynomials \cite{KostantKumar1}, \cite{KostantKumar2}, \cite{Kumar}, \cite{Billey}. In \cite{Kumar} S. Kumar showed that if $w_1$ and~$w_2$ are arbitrary elements of~$W$ and $d_{w_1}\neq d_{w_2}$, then $C_{w_1}\neq C_{w_2}$. We give three equivalent definitions of Kostant--Kumar polynomials and formulate their properties needed for the sequel. In Section~\ref{sect:divisibility}, we recall basic definitions and facts about root systems of type $E$ and prove the main technical fact about divisibility of Kostant--Kumar polynomials, see Proposition~\ref{prop:crucial}. Finally, Section~\ref{sect:good_pairs} contains the notion of a good pair of involutions and the proof of our main result, Theorem~\ref{mtheo:non_red}, based on Proposition~\ref{prop:crucial} and detailed consideration of configurations of roots, see Proposition~\ref{prop:good_pair}.

\medskip\textsc{Acknowledgements}. the Foundation for the Advancement of Theoretical Physics and Mathe\-matics ``BASIS'', grant no. 18--1--7--2--1. The first author was also partially supported by RFBR grant no. 20--01--00091a. These foundations are gratefully acknowledged.

\sect{Kostant--Kumar polynomials}\label{sect:Kostant_Kumar}

\sst Let $w$ be an element of the Weyl group $W$. Here we recall the precise definition of the Kostant--Kumar polynomial $d_w$, explain how to compute it in combinatorial terms, and show that it depends only on the scheme structure of $C_w$, see~\cite{Kumar} for the details.

The torus $T$ acts on the Schubert variety $X_w$ by left multiplications (or, equivalently, by con\-ju\-ga\-tions). The point $p$ is invariant under this action, hence there is the structure of a $T$-module on the local ring $\Ou$. The action of $T$ on $\Ou$ preserves the filtration by powers of the ideal $\mt$, so we obtain the structure of a $T$-module on the algebra $R=\mathrm{gr}\,\Ou$. By~\cite[Theorem~2.2]{Kumar}, $R$ can be decomposed into a direct sum of its finite-dimensional weight subspaces: $$R=\bigoplus_{\lambda\in\Xt(T)}R_{\lambda}.$$
Here $\htt$ is the Lie algebra of the torus $T$, $\Xt(T)\subseteq\htt^*$ is the character lattice of $T$ and $R_{\lambda}=\{f\in R\mid t.f=\lambda(t)f\}$ is the weight subspace of weight $\lambda$. Let $\Lambda$ be the $\Zp$-module consisting of all (possibly infinite) $\Zp$-linear combinations of linearly independent elements $e^{\lambda}$, $\lambda\in\Xt(T)$. The \emph{formal character} of~$R$ is an element of $\Lambda$ of the form $$\chr R=\sum_{\lambda\in\Xt(T)}m_{\lambda}e^{\lambda},$$ where $m_{\lambda}=\dim R_{\lambda}$.\newpage

Now, pick an element $a=\sum_{\lambda\in\Xt(T)}n_{\lambda}e^{\lambda}\in\Lambda$. Assume that there are finitely many $\lambda\in\Xt(T)$ such that $n_{\lambda}\neq0$. Given $k\geq0$, one can define the polynomial $$[a]_k=\sum_{\lambda\in\Xt(T)}n_{\lambda}\cdot\dfrac{\lambda^k}{k!}\in S=\Cp[\htt].$$ Denote $[a]=[a]_{k_0}$, where $k_0$ is minimal among all non-negative numbers $k$ such that $[a]_k\neq0$. For instance, if $a=1-e^{\lambda}$, then $[a]_0=0$ and $[a]=[a]_1=-\lambda$ (here we denote $1=e^0$).

Let $A$ be the submodule of $\Lambda$ consisting of all finite linear combinations. It is a commutative ring with respect to the multiplication $e^{\lambda}\cdot e^{\mu}=e^{\lambda+\mu}$. In fact, it is just the group ring of $\Xt(T)$. Denote the field of fractions of the ring $A$ by $Q$. To each element of $Q$ of the form $q=a/b$, $a$,~$b\in A$, one can assign the element $$[q]=\dfrac{[a]}{[b]}\in\Cp(\htt)$$ of the field of rational functions on $\htt$. Note that this element is well-defined~\cite{Kumar}.

There exists an involution $q\mapsto q^*$ on $Q$ defined by $$e^{\lambda}\mapsto (e^{\lambda})^*=e^{-\lambda}.$$ It turns out \cite[Theorem 2.2]{Kumar} that the character $\chr R$ belongs to $Q$, hence $(\chr R)^*\in Q$, too. (One can consider the field $Q$ of rational functions as a subring of the ring $\Lambda$.) Finally, we put $$c_w=[(\chr R)^*],\ d_w=(-1)^{l(w)}\cdot c_w\cdot\prod_{\alpha\in\Phi^+}\alpha.$$ Here $l(w)$ is the length of $w$ in the Weyl group $W$ with respect to the set of simple roots $\Delta$. Evidently, $c_w$ and $d_w$ belong to $\Cp(\htt)$; in fact, $d_w$ is a polynomial, i.e., it belongs to the algebra $S=\Cp[\htt]$ of regular functions on $\htt$, see \cite{KostantKumar2} and \cite[Theorem 7.2.6]{BilleyLakshmibai}. \defi{Let $w$ be an element of the Weyl group $W$. The polynomial $d_w\in S$ is called the \emph{Kostant--Kumar polynomial} associated with $w$.}

It follows from the definition that $c_w$ and $d_w$ depend only on the canonical structure of a $T$-module on the algebra $R$ of regular functions on the tangent cone $C_w$. Thus, to prove that the tangent cones corresponding to elements $w_1$, $w_2$ of the Weyl group are distinct, it is enough to check that ${c_{w_1}\neq c_{w_2}}$, or, equivalently, $d_{w_1}\neq d_{w_2}$.

On the other hand, there is a purely combinatorial description of Kostant--Kumar polynomials. To give this description, we need some more notation. Let $w$, $v$ be elements of~$W$. Fix a reduced decomposition of the element $w=s_{i_1}\ldots s_{i_l}$. (Here $\alpha_1,\ldots,\alpha_n\in\Delta$ are simple roots and $s_i=s_{\alpha_i}$ is the simple reflection corresponding to $\alpha_i$.) Put
\begin{equation*}
c_{w,v}=(-1)^{l(w)}\cdot\sum\dfrac{1}{s_{i_1}^{\epsi_1}\alpha_{i_1}}\cdot\dfrac{1}{s_{i_1}^{\epsi_1}
s_{i_2}^{\epsi_2}\alpha_{i_2}}\cdot\ldots\cdot\dfrac{1}{s_{i_1}^{\epsi_1}\ldots s_{i_l}^{\epsi_l}\alpha_{i_l}},
\end{equation*}
where the sum is taken over all sequences $(\epsi_1,\ldots,\epsi_l)$ of zeroes and units such that ${s_{i_1}^{\epsi_1}\ldots s_{i_l}^{\epsi_l}=v}$. Actually, the element $c_{w,v}\in\Cp(\htt)$ depends only on $w$ and $v$, not on the choice of a reduced de\-com\-po\-si\-tion of $w$ \cite[Section 3]{Kumar}.

\exam{Let $\Phi=A_n$. Put $w=s_1s_2s_1$. To compute $c_{w,\id}$, we should take the sum over two sequences, $(0,0,0)$ and $(1,0,1)$. Hence
\begin{equation*}
c_{w,\id}=(-1)^3\cdot\left(\dfrac{1}{\alpha_1\alpha_2\alpha_1}+\dfrac{1}{-\alpha_1(\alpha_1+\alpha_2)\alpha_1}\right)
=-\dfrac{1}{\alpha_1\alpha_2(\alpha_1+\alpha_2)}.
\end{equation*}}

A remarkable fact is that $c_{w,\id}=c_w$, hence to prove that the tangent cones to Schubert varieties do not coincide as subschemes, we need only combinatorics of the Weyl group. Note also that for classical Weyl groups, elements $c_{w,v}$ are closely related to Schubert polynomials \cite{Billey}.

Finally, we will present an original definition of elements $c_{w,v}$ using so-called nil-Hecke ring (see \cite{Kumar} and~\cite[Section 7.1]{BilleyLakshmibai}). The group $W$ naturally acts on $\Cp(\htt)$ by automorphisms. Denote by $Q_W$ the vector space over $\Cp(\htt)$ with basis $\{\delta_w,\ w\in W\}$. It is a ring with respect to the multiplication $$f\delta_v\cdot g\delta_w=fv(g)\delta_{vw}.$$
This ring is called the \emph{nil-Hecke ring}. To each $i$ from 1 to $n$ put $$x_i=\alpha_i^{-1}(\delta_{s_i}-\delta_{\id}).$$
Let $w\in W$ and $w=s_{i_1}\ldots s_{i_l}$ be a reduced decomposition of $w$. Then the element $$x_w=x_{i_1}\ldots x_{i_l}$$ does not depend on the choice of a reduced decomposition of $w$ \cite[Proposition 2.1]{KostantKumar1}.

Moreover, it turns out that $\{x_w,\ w\in W\}$ is a $\Cp(\htt)$-basis of~$Q_W$ \cite[Proposition 2.2]{KostantKumar1}, and
\begin{equation*}
x_w=\sum\nolimits_{v\in W}c_{w,v}\delta_v.\\
\end{equation*}
Actually, if $w,v\in W$, then
\begin{equation}
\begin{split}
&\text{a) }x_v\cdot x_w=\begin{cases}x_{vw},&\text{ if }l(vw)=l(v)+l(w),\\
0,&\text{ otherwise},
\end{cases}\\
&\text{b) }c_{w,v}=-v(\alpha_i)^{-1}(c_{ws_i,v}+c_{ws_i,vs_i}),\text{ if }l(ws_i)=l(w)-1,\\
&\text{c) }c_{w,v}=\alpha_i^{-1}(s_i(c_{s_iw,s_iv})-c_{s_iw,v}),\text{ if }l(s_iw)=l(w)-1.\\
\end{split}\label{formula_x_v_x_w}
\end{equation}
The first property is proved in \cite[Proposition 2.2]{KostantKumar1}. The second and the third properties follow immediately from the first one and the definitions (see also the proof of \cite[Corollary 3.2]{Kumar}).

\nota{Suppose~$\Phi$ is a union of its\label{nota:irred} subsystems $\Phi_1$ and $\Phi_2$ contained in mutually orthogonal subspaces. Let $W_1$, $W_2$ be the Weyl groups of $\Phi_1$, $\Phi_2$ respectively, so $W=W_1\times W_2$. Denote $\Delta_1=\Delta\cap\Phi_1=\{\alpha_1,\ldots,\alpha_r\}$ and $\Delta_2=\Delta\cap\Phi_2=\{\beta_1,\ldots,\beta_s\}$, then $$\Cp[\htt]\cong\Cp[\alpha_1,\ldots,\alpha_r,\beta_1,\ldots,\beta_s].$$
Given $v\in W_i$, $i=1,~2$, denote by $d_v^i$ its Kostant--Kumar polynomial. We can consider $d_v^1$ (re\-spec\-tively,~$d_v^2$) as an element of~$\Cp[\htt]$ depending only on $\alpha_1,\ldots,\alpha_r$ (respectively, on $\beta_1,~\ldots,~\beta_s$). We define $c_v^i\in\Cp(\htt)$, $i=1,~2$, by a similar way. Let $w\in W$, $w_1\in W_1$, $w_2\in W_2$ and $w=w_1w_2$. Repeating literally the proof of \cite[Proposition 1.6]{EliseevIgnatyev}, we obtain the following: $$d_w=d_{w_1}^1d_{w_2}^2,\ c_w=c_{w_1}^1c_{w_2}^2.$$ Thus, to prove Theorem~\ref{mtheo:non_red} it is enough to prove this theorem for irreducible root systems of type $E$, because $\Cp[\htt]$ is a unique factorization domain.}

\newpage\sect{Divisibility in $\Cp[\htt]$}\label{sect:divisibility}

Throughout this section, $\Phi$ denotes an irreducible root system of type $E_6$, $E_7$ or $E_8$. Below we briefly recall some facts about $\Phi$. (We follow the notation from \cite{Bourbaki}.) Let $\epsi_1$, $\ldots$, $\epsi_n$ be the standard basis of the Euclidean space~$\Rp^n$. As usual, we identify the set $\Phi^+$ of positive roots with the following subset of $\Rp^n$:
\begin{equation*}
\begin{split}
E_6^+&=\{(\pm\epsilon_i + \epsilon_j),~ 1\leq i < j \leq 5 \} \cup \\ & \cup \left\{  \frac{1}{2} \left( \epsilon_8-\epsilon_7-\epsilon_6+\sum \limits_{i=1}^{5} (-1)^{\nu(i)} \epsilon_i \right),~\sum \limits_{i=1}^{5} \nu(i) \text{ is even} \right\},\\
E_7^+&=\{(\pm\epsilon_i+ \epsilon_j) 1\leq i < j \leq 6 \} \cup \{ (\epsilon_7-\epsilon_8) \} \cup \\ & \cup \left\{ \frac{1}{2} \left( \epsilon_7-\epsilon_8+\sum \limits_{i=1}^{6} (-1)^{\nu(i)} \epsilon_i \right),~\sum \limits_{i=1}^{6} \nu(i) \text{ is even} \right\},\\
E_8^+&=\{\pm\epsilon_i+ \epsilon_j 1\leq i < j \leq 8 \} \cup \left\{ \frac{1}{2} \sum \limits_{i=1}^{8} (-1)^{\nu(i)} \epsilon_i,~\sum \limits_{i=1}^{8} \nu(i) \text{ is even} \right\},\\
\end{split}
\end{equation*}
so $W$ can be considered as a subgroup of the orthogonal group $O(\Rp^n)$.

The simple roots have the following form.
\begin{equation}\label{formula:simple_roots}
\begin{split}
\Phi=E_6\colon\alpha_1&=\frac{1}{2}(\epsilon_1+\epsilon_8)-\frac{1}{2}(\epsilon_2+\epsilon_3+\epsilon_4+\epsilon_5+\epsilon_6+\epsilon_7),~\alpha_2=\epsilon_1+\epsilon_2,~\alpha_3=\epsilon_2-\epsilon_1,\\
\alpha_4&=\epsilon_3-\epsilon_2,~\alpha_5=\epsilon_4-\epsilon_3,~\alpha_6=\epsilon_5-\epsilon_4;\\
\Phi=E_7\colon\alpha_1&=\frac{1}{2}(\epsilon_1+\epsilon_8)-\frac{1}{2}(\epsilon_2+\epsilon_3+\epsilon_4+\epsilon_5+\epsilon_6+\epsilon_7),~\alpha_2=\epsilon_1+\epsilon_2,~\alpha_3=\epsilon_2-\epsilon_1,\\
\alpha_4&=\epsilon_3-\epsilon_2,~\alpha_5=\epsilon_4-\epsilon_3,~\alpha_6=\epsilon_5-\epsilon_4,\alpha_7=\epsilon_6-\epsilon_5.\\
\Phi=E_8\colon\alpha_1&=\frac{1}{2}(\epsilon_1+\epsilon_8)-\frac{1}{2}(\epsilon_2+\epsilon_3+\epsilon_4+\epsilon_5+\epsilon_6+\epsilon_7),~\alpha_2=\epsilon_1+\epsilon_2,~\alpha_3=\epsilon_2-\epsilon_1,\\
\alpha_4&=\epsilon_3-\epsilon_2,~\alpha_5=\epsilon_4-\epsilon_3,~\alpha_6=\epsilon_5-\epsilon_4,\alpha_7=\epsilon_6-\epsilon_5,~\alpha_8=\epsilon_7-\epsilon_6.\\
\end{split}
\end{equation}

We say that $v$ is less or equal to $w$ with respect to the \emph{Bruhat order}, written $v\leq w$, if some reduced decomposition for $v$ is a subword of some reduced decomposition for $w$. It is well-known that this order plays the crucial role in many geometric aspects of theory of algebraic groups. For instance, the Bruhat order encodes the incidences among Schubert varieties, i.e., $X_v$ is contained in $X_w$ if and only if $v\leq w$. It turns out that $c_{w,v}$ is non-zero if and only if $v\leq w$ \cite[Corollary 3.2]{Kumar}. For example, $c_w=c_{w,\id}$ is non-zero for \emph{any} $w$, because $\id$ is the smallest element of $W$ with respect to the Bruhat order. Note that given $v,w\in W$, there exists $g_{w,v}\in S=\Cp[\htt]$ such that
\begin{equation}
c_{w,v}=g_{w,v}\cdot\prod_{\alpha>0,~s_{\alpha}v\leq w}\alpha^{-1},\label{formula:dyer}
\end{equation}
see \cite{Dyer} and \cite[Theorem 7.1.11]{BilleyLakshmibai}

Since we fixed the order on the set of simple roots, one can consider the lexicographic total order on the set of positive roots: given $\alpha=\sum a_i\alpha_i$ and $\beta=\sum b_i\alpha_i$, we write $\alpha\prec\beta$ if there exists $j$ such that $a_i=b_i$ for all $i<j$ and $a_i<b_i$. Let $w$ be an involution in the Weyl group $W$ of $\Phi$. Denote by $s_{\alpha}$ the reflection in $W$ corresponding to a root $\alpha$. Denote by $\beta_1$ the maximal (with respect to the order~$\preceq$) root among all roots $\beta\in\Phi^+$ for which $w(\beta)=-\beta$. Next, for $i\geq1$, denote by $\beta_{i+1}$ the maximal root among all roots $\beta\in\Phi^+$ such that $w_i(\beta)=-\beta$, where $$w_i=s_{\beta_i}\circ s_{\beta_{i-1}}\circ\ldots\circ s_{\beta_1}\circ w.$$ One can easily check that $w_k$ coincides with the identity element of $W$ for certain $k$.\newpage

\defi{The set $\Supp{w}=\{\beta_1,~\ldots,~\beta_k\}$ is called the \emph{support} $\Supp{\sigma}$ of $w$. It turns out that $\Supp{w}$ is an orthogonal subset of $\Phi^+$ \cite[Theorem 5.4]{Springer}. Note that $$w=\prod\nolimits_{\beta\in\Supp{w}}s_{\beta},$$ where the product is taken in any fixed order.}

\lemmp{Let $w_1$, $w_2$ be\label{lemm:Supp} involutions in $W$. If $\Supp{w_1}\subset\Supp{w_2}$ then $w_1\leq w_2$.}{The well-known Strong Exchange Condition (see, e.g., \cite[Proposition 3.1 (ii)]{Deodhar}) implies that, given $w\in W$ and $\alpha\in\Phi$, one has $l(ws_{\alpha})>l(w)$ if and only if $w{\alpha}\in\Phi^+$. On the other hand (see, e.g., \cite[Definition 2.1.1]{BjornerBrenti}), $l(ws_{\alpha})>l(w)$ if and only if $ws_{\alpha}>w$. Hence, $w(\alpha)\in\Phi^+$ if and only if $ws_{\alpha}>w$. Let $\Supp{w_2}\setminus\Supp{w_1}=\{\beta_1,~\ldots,~\beta_k\}$, then
\begin{equation*}
w_2=w_1\cdot\prod_{\beta\in\Supp{w_2}\setminus\Supp{w_1}}s_{\beta}=w_1s_{\beta_1}\ldots s_{\beta_k}.
\end{equation*}
Next, denote $v_i=w_1s_{\beta_1}\ldots s_{\beta_i-1}$ for $1\leq i\leq k+1$, so that $v_1=w_1$ and $v_{k+1}=w_2$. Then, clearly, $v_i(\beta_{i})=\beta_i\in\Phi^+$, thus, $w_1=v_1<v_2<\ldots<v_{k+1}=w_2$, as required.}

\defi{The subset $\Co_1=\{\beta\in\Phi^+\mid\alpha_1\preceq\beta\}$ is called the \emph{first column of} $\Phi^+$.}

%Set $\beta_0$ to be the highest root in $\Phi^+$ and denote
%\begin{equation*}
%\Co=\begin{cases}\Co_1\text{ or }\Co_6,&\text{if }\Phi=E_6,\\
%\Co_7,&\text{if }\Phi=E_7,\\
%\Co_8\setminus\{\beta_0\},&\text{if }\Phi=E_8.\\
%\end{cases}
%\end{equation*}

We will essentially use the following standard fact about parabolic subgroups of the Weil group $W$.

\mtheo{\textup{\cite[Proposition 1.10 (c)]{Humpreys2}} Let $I$ be a subset\label{theo:parabolic} of the set $\Delta$ of simple roots. Denote by $W_I$ the parabolic subgroup of $W$ generated by the simple reflections $s_{\alpha}$\textup, $\alpha\in I$. Put also\break $W^I=\{w\in W\mid l(ws_{\alpha})>l(w)\text{ for all }\alpha\in I\}$. Given $w\in W$\textup, there exist unique $u\in W^I$ and $v\in W_I$ such that $w=uv$. Their lengths satisfy $l(w)=l(u)+l(v)$.}

The following proposition plays the crucial role in the proof of the main result (cf. \cite[Lem\-mas~2.4,~2.5]{EliseevIgnatyev}, \cite[Lemma 2.6]{BochkarevIgnatyevShevchenko} and \cite[Lemma 2.7]{IgnatyevShevchenko1}).

\propp{Let $w\in W$ be an involution. Assume that $\Supp{w}\cap\Co_1=\{\beta\}$ and the reflection $s_{\beta}$ has a reduced decomposition of the form $s_{\beta}=u_{\beta}v_{\beta}$ for a certain element $v_{\beta}$ from the subgroup $\wt W$ of~$W$ generated by the reflections $s_i$\textup, $i\neq1$\textup, so that $u_{\beta}=v_{\beta}^{-1}s_1$ and $l(u_{\beta}s_i)=l(u)+1$ for all $i\neq1$. Then $\beta$ does not divide $d_w$ in $\Cp[\htt]$.\label{prop:crucial}}{%Denote by $\wt\Phi$ the root system corresponding to $\wt W$; in fact, $\wt\Phi^+=\Phi^+\setminus\Co_r$. Denote by $\wt d_w\in\wt S=\Cp[\alpha_i,~i\neq r]$ the Kostant--Kumar polynomial of $w$ considered as an element of~$\wt W$; define $\wt c_w\in\Cp(\alpha_i,~i\neq r)$ by the similar way. Since $\wt W$ is a parabolic subgroup of $W$, the length of $w$ as an element of $\wt W$ equals the length of $w$ as an element of $W$. Further, any reduced decomposition for $w$ in $\wt W$ is a reduced decomposition for $w$ in $W$. This means that $\wt c_w=c_w$, so
Denote
\begin{equation*}
\wt W^1=\{w\in W\mid l(ws_i)=l(w)+1\text{ for all }i\neq1\}=\{w\in W\mid w(\alpha_i)\in\Phi^+\text{ for all }i\neq1\}.
\end{equation*}
Applying Theorem~\ref{theo:parabolic} to the subset $I=\Delta\setminus\{\alpha_1\}$, we see that there exist unique $u\in\wt W^1=W^I$ and $v\in\wt W=W_I$ such that $w=uv$. We claim that in fact $u=u_{\beta}$. Indeed, denote $w'=\prod_{\alpha\in\Supp{w},~\alpha\neq\beta}s_{\alpha}$, then one can write
\begin{equation*}
w=\prod_{\alpha\in\Supp{w}}s_{\alpha}=s_{\beta}w'=u_{\beta}v_{\beta}w'.
\end{equation*}
But $v_{\beta}w'\in\wt W$, while $u_{\beta}\in\wt W^1$, which means that $u=u_{\beta}$ and $v=v_{\beta}w'$.

We claim that
\begin{equation}
\begin{split}
c_w&=-\dfrac{c_{us_1,g_0}g_0(c_{v,g_0^{-1}})}{\beta}-\sum_{g\leq u,\ g^{-1}\leq v,\ g\neq g_0}\dfrac{c_{us_1,g}g(c_{v,g^{-1}})}{g(\alpha_1)}\\
&=\beta^{-1}\cdot g_0(c_{v,g_0^{-1}})\cdot\dfrac{K}{L}+\dfrac{M}{N}\label{formula:KLMN}
\end{split}
\end{equation}
Here $g_0=us_1$ and $K,L$ and $M,N\in\Cp[\htt]$ are pairs of coprime polynomials such that the root $\beta$ (considered as an element of $\Cp[\htt]$) divides neither $K$ nor $N$.

Indeed, one can prove (\ref{formula:KLMN}) using (\ref{formula_x_v_x_w}) and arguing as in the proof of
\cite[Lemma 2.5]{EliseevIgnatyev}. Namely, since $l(w)=l(u)+l(v)$, formula (\ref{formula_x_v_x_w}a) shows that
\begin{equation*}
\begin{split}
x_w&=\sum\nolimits_{s\in W}c_{w,s}\delta_s=x_ux_v=\sum\nolimits_{g,h\in W}c_{u,g}\delta_g\cdot c_{v,h}\delta_h\\
&=\sum\nolimits_{g,h\in W}c_{u,g}g(c_{v,h})\delta_{gh}=\sum\nolimits_{s\in W}\left(\sum\nolimits_{g\in W}c_{u,g}g(c_{v,g^{-1}s})\right)\delta_s.
\end{split}
\end{equation*}
Thus, for any $s\in W$, the coefficient of $\delta_s$ is equal to
$$c_{w,s}=\sum_{g\in W}c_{u,g}g(c_{v,g^{-1}s}),$$
in particular,
$$c_w=c_{w,\id}=\sum_{g\in W}c_{u,g}g(c_{v,g^{-1}}).$$
Moreover, since $c_{p,q}\neq0$ if and only if $p\geq q$, the sum in the right hand side is taken over permutations $g$ such that $u\geq g$ and $v\geq g^{-1}$. Denote the set of such permutations by $U$. Note that $g\in U$ implies that $g$ is obtained from $u=v_{\beta}^{-1}s_1$ by deleting $s_1$ and, possibly, some other simple reflections. (If $s_1$ is not deleted, then the condition $v\geq g^{-1}$ does not hold.) Hence
$$c_w=c_{w,\id}=\sum_{g\in U}c_{u,g}g(c_{v,g^{-1}}).$$

Using (\ref{formula_x_v_x_w}b) and the fact that $l(us_1)=l(u)-1$, we obtain
$$c_{u,g}=-g(\alpha_1)^{-1}(c_{us_1,g}+c_{us_1,gs_1})=-g(\alpha_1)^{-1}c_{us_1,g},$$
because $us_1\not\geq gs_1$ and so $c_{us_1,gs_1}=0$. Thus,
$$c_w=-\sum_{g\in U}\dfrac{c_{us_1,g}g(c_{v,g^{-1}})}{g(\alpha_1)}.$$
It is easy to check that there is at most one element $g$ such that $g(\alpha_1)=\beta$ and $g\in U$, namely, the element $g_0=us_1=v_{\beta}^{-1}$. Indeed, $$s_{\beta}=v_{\beta}^{-1}s_1v_{\beta}=s_{v_{\beta}^{-1}(\alpha_1)},$$ hence $v_{\beta}^{-1}(\alpha_1)=\pm\beta$. But $v_{\beta}^{-1}$ belongs to~$\wt W$, consequently, $v_{\beta}^{-1}(\alpha_1)=\alpha_1+\ldots\in\Phi^+$. We conclude that $g_0=v_{\beta}^{-1}$ sends $\alpha_1$ to $\beta$. On the other hand, if $g(\alpha_1)=\beta$ for some $g\neq v_{\beta}^{-1}$ from $U$, then $s_{\beta}=v_{\beta}^{-1}s_1v_{\beta}$ is not a reduced decomposition of $s_{\beta}$, a contradiction.

Assume for a moment that $g_0$ belongs to $U$, i.e., $v\geq g_0^{-1}$. Then
\begin{equation}
c_w=-\dfrac{c_{us_1,g_0}g_0(c_{v,g_0^{-1}})}{\beta}-\sum_{g\in U,\ g\neq g_0}\dfrac{c_{us_1,g}g(c_{v,g^{-1}})}{g(\alpha_1)}.\label{formula:c_w_c_v}
\end{equation} By $S'$ (resp. $Q'$) denote the subalgebra of $S=\Cp[\htt]$ (resp. the subfield of $\Cp(\htt)$) generated by\linebreak $\alpha_i$, $i\neq1$, then $c_{v,g_0^{-1}}\in Q'$ because $v$, $g_0^{-1}\in\wt W$. Since $g\in\wt W$, $g(c_{v,g_0^{-1}})\in Q'$, too. In particular, if $g(c_{v,g_0^{-1}})=G_1/G_2$ and $G_1$, $G_2\in S'$ are coprime, then $\beta$ does not divide $G_1$. On the other hand, $c_{us_1,g_0}\in Q'$, because both $us_1=g_0$ belongs to $\wt W$. We conclude that the first summand in (\ref{formula:c_w_c_v}) has the form $g_0(c_{v,g_0^{-1}})\cdot K/\beta L$ for some coprime $K,L\in\Cp[\htt]$. Finally, if $g\in U$ and $g\neq g_0$, then $g(c_{v,g^{-1}})\in Q'$. Again, since $us_1$ and $g$ belong to $\wt W$, one has $c_{us_1,g}\in Q'$. We see that if the latter sum in~(\ref{formula:c_w_c_v}) is equal to $M/N$, where $M,N\in\Cp[\htt]$ are coprime, then $\beta$ does not divide $N$.

To prove that $\beta$ does not divide $d_w$, it is enough to show that $c_{v,g_0^{-1}}\neq0$, i.e., $v\geq g_0^{-1}$ (or, equivalently, $v^{-1}\geq g_0$). By Lemma~\ref{lemm:Supp}, $s_{\beta}\leq w$, According to \cite[Chapter 2, Exercise 21]{BjornerBrenti}, this is equivalent to $v_{\beta}\leq v$. Hence, $g_0=v_{\beta}^{-1}\leq v^{-1}$, which concludes the proof.}

\newpage\sect{Good pairs\label{sect:good_pairs} of involutions} In   this section, we formulate and prove the main result of the paper, Theorem~\ref{mtheo:non_red}. To do this, we need to introduce the notion of a good pair of involutions. Recall the set of simple roots from (\ref{formula:simple_roots}). We will order the simple roots as follows.
\begin{center}
\begin{tabular}{|l|l|}
\hline
Type of $\Phi$&Order of simple roots\\
\hline\hline
$E_6$&$\alpha_1,~\alpha_2,~\alpha_3,~\alpha_4,~\alpha_5,~\alpha_6$ or\\
&$\alpha_2,~\alpha_6,~\alpha_3,~\alpha_5,~\alpha_4,~\alpha_1$\\
\hline
$E_7$&$\alpha_3,~\alpha_7,~\alpha_4,~\alpha_6,~\alpha_5,~\alpha_2,~\alpha_1$\\
\hline
$E_8$&$\alpha_4,~\alpha_8,~\alpha_5,~\alpha_7,~\alpha_6,~\alpha_3,~\alpha_2,~\alpha_1$\\
\hline
\end{tabular}
\end{center}
Since the set of simple roots is ordered, the support of an involution and the first column are well-defined.

\defi{Let $w_1$, $w_2$ be\label{defi:good_pair} involutions in $W$. We say that they form a \emph{good pair of involutions} if $\Supp{w_i}\cap\Co_1=\{\beta_i\}$ for $i=1,~2$ such that $\beta_1\neq\beta_2$, both $\beta_1$ and $\beta_2$ are not maximal in $C_1$ for $\Phi=E_8$, and $s_{\beta_1}\nleq w_2$ or $s_{\beta_2}\nleq w_1$.}

Now we are ready to prove our main result, Theorem~\ref{mtheo:non_red}, which claims that if $w_1$, $w_2$ is a good pair of involutions then the corresponding tangent cones $C_{w_1}$ and $C_{w_2}$ do not coincide. This follows immediately from the following proposition.

\propp{Let $w_1$, $w_2$ be a good\label{prop:good_pair} pair of involutions in $W$. Then $d_{w_1}\neq d_{w_2}$.}{Let $\beta=\beta_1$ or $\beta_2$, and $s_{\beta}=uv_{\beta}$ be as in Proposition~\ref{prop:crucial}. In the tables below we list the elements $u$ for all possible $\beta$. The first (respectively, the second) column of the table contains the sequence $(c_1,~\ldots,~c_8)$ (respectively, $(b_1,~\ldots,~b_n)$) if $$\beta=\sum_{i=1}^8c_i\epsi_i=\sum_{i=1}^n\alpha_i,~n=\rk\Phi.$$ The third column contains a reduced decomposition of $u$.
	\begin{longtable}{|c|c|c|}
		\hline
\multicolumn{3}{|c|}{Case $\Phi=E_6$ with the order $\alpha_1,~\alpha_2,~\alpha_3,~\alpha_4,~\alpha_5,~\alpha_6$}\\\hline
		$(c_1,~\ldots,~c_8)$ & $(b_1,~\ldots,~b_6)$& Reduced decomposition of $u$\\
		\hline\hline
		$(\frac{1}{2}, -\frac{1}{2}, -\frac{1}{2}, -\frac{1}{2}, -\frac{1}{2}, -\frac{1}{2}, -\frac{1}{2}, \frac{1}{2})$ & $100000$ & $s_1$  \\
		\hline
		$(-\frac{1}{2}, \frac{1}{2}, -\frac{1}{2}, -\frac{1}{2}, -\frac{1}{2}, -\frac{1}{2}, -\frac{1}{2}, \frac{1}{2})$ & $101000$ & $s_3s_1$\\
		\hline
		$(-\frac{1}{2}, -\frac{1}{2}, \frac{1}{2}, -\frac{1}{2}, -\frac{1}{2}, -\frac{1}{2}, -\frac{1}{2}, \frac{1}{2})$ & $101100$ & $s_4s_3s_1$\\
		\hline
		$(\frac{1}{2}, \frac{1}{2}, \frac{1}{2}, -\frac{1}{2}, -\frac{1}{2}, -\frac{1}{2}, -\frac{1}{2}, \frac{1}{2})$ & $111100$ & $s_2s_4s_3s_1$\\
		\hline
		$(-\frac{1}{2}, -\frac{1}{2}, -\frac{1}{2}, \frac{1}{2}, -\frac{1}{2}, -\frac{1}{2}, -\frac{1}{2}, \frac{1}{2})$ & $101110$ & $s_5s_4s_3s_1$ \\
		\hline
		$(-\frac{1}{2}, -\frac{1}{2}, -\frac{1}{2}, -\frac{1}{2}, \frac{1}{2}, -\frac{1}{2}, -\frac{1}{2}, \frac{1}{2})$ & $101111$ & $s_6s_5s_4s_3s_1$ \\
		\hline
		$(\frac{1}{2}, \frac{1}{2}, -\frac{1}{2}, \frac{1}{2}, -\frac{1}{2}, -\frac{1}{2}, -\frac{1}{2}, \frac{1}{2})$ & $111110$ & $s_5s_2s_4s_3s_1$ \\
		\hline
		$(\frac{1}{2}, -\frac{1}{2}, \frac{1}{2}, \frac{1}{2}, -\frac{1}{2}, -\frac{1}{2}, -\frac{1}{2}, \frac{1}{2})$ & $111210$ & $s_4s_5s_2s_4s_3s_1$ \\
		\hline
		$(\frac{1}{2}, \frac{1}{2}, -\frac{1}{2}, -\frac{1}{2}, \frac{1}{2}, -\frac{1}{2}, -\frac{1}{2}, \frac{1}{2})$ & $111111$ & $s_6s_5s_2s_4s_3s_1$ \\
		\hline
		$(-\frac{1}{2}, \frac{1}{2}, \frac{1}{2}, \frac{1}{2}, -\frac{1}{2}, -\frac{1}{2}, -\frac{1}{2}, \frac{1}{2})$ & $112210$ & $s_3s_4s_5s_2s_4s_3s_1$  \\
		\hline
		$(\frac{1}{2}, -\frac{1}{2}, \frac{1}{2}, -\frac{1}{2}, \frac{1}{2}, -\frac{1}{2}, -\frac{1}{2}, \frac{1}{2})$ & $111211$ & $s_6s_4s_5s_2s_4s_3s_1$ \\
		\hline
		$(-\frac{1}{2}, \frac{1}{2}, \frac{1}{2}, -\frac{1}{2}, \frac{1}{2}, -\frac{1}{2}, -\frac{1}{2}, \frac{1}{2})$ & $112211$ & $s_3s_6s_4s_5s_2s_4s_3s_1$  \\
		\hline
		$(\frac{1}{2}, -\frac{1}{2}, -\frac{1}{2}, \frac{1}{2}, \frac{1}{2}, -\frac{1}{2}, -\frac{1}{2}, \frac{1}{2})$ & $111221$ & $s_5s_6s_4s_5s_2s_4s_3s_1$  \\
		\hline
		$(-\frac{1}{2}, \frac{1}{2}, -\frac{1}{2}, \frac{1}{2}, \frac{1}{2}, -\frac{1}{2}, -\frac{1}{2}, \frac{1}{2})$ & $112221$ & $s_3s_5s_6s_4s_5s_2s_4s_3s_1$  \\
		\hline
		$(-\frac{1}{2}, -\frac{1}{2}, \frac{1}{2}, \frac{1}{2}, \frac{1}{2}, -\frac{1}{2}, -\frac{1}{2}, \frac{1}{2})$ & $112321$ & $s_4s_3s_5s_6s_4s_5s_2s_4s_3s_1$ \\
		\hline
		$(\frac{1}{2}, \frac{1}{2}, \frac{1}{2}, \frac{1}{2}, \frac{1}{2}, -\frac{1}{2}, -\frac{1}{2}, \frac{1}{2})$ & $122321$ & $s_2s_4s_3s_5s_6s_4s_5s_2s_4s_3s_1$  \\
		\hline
	\end{longtable}

\begin{longtable}{|c|c|c|}
		\hline
\multicolumn{3}{|c|}{Case $\Phi=E_6$ with the order $\alpha_2,~\alpha_6,~\alpha_3,~\alpha_5,~\alpha_4,~\alpha_1$}\\\hline
		$(c_1,~\ldots,~c_8)$ & $(b_1,~\ldots,~b_6)$& Reduced decomposition of $u$\\
		\hline\hline
		$(0, 0, 0, -1, 1, 0, 0, 0)$ & $000001$ & $s_6$ \\
		\hline
		$(0, 0, -1, 0, 1, 0, 0, 0)$ & $000011$ & $s_5s_6$ \\
		\hline
		$(0, -1, 0, 0, 1, 0, 0, 0)$ & $000111$ & $s_4s_5s_6$ \\
		\hline
		$(1, 0, 0, 0, 1, 0, 0, 0)$ & $010111$ & $s_2s_4s_5s_6$ \\
		\hline
		$(-1, 0, 0, 0, 1, 0, 0, 0)$ & $001111$ & $s_3s_4s_5s_6$ \\
		\hline
		$(-\frac{1}{2}, -\frac{1}{2}, -\frac{1}{2}, -\frac{1}{2}, \frac{1}{2}, -\frac{1}{2}, -\frac{1}{2}, \frac{1}{2})$ & $101111$ & $s_1s_3s_4s_5s_6$ \\
		\hline
		$(0, 1, 0, 0, 1, 0, 0, 0)$ & $011111$ & $s_2s_3s_4s_5s_6$ \\
		\hline
		$(0, 0, 1, 0, 1, 0, 0, 0)$ & $011211$ & $s_4s_2s_3s_4s_5s_6$ \\
		\hline
		$(\frac{1}{2}, \frac{1}{2}, -\frac{1}{2}, -\frac{1}{2}, \frac{1}{2}, -\frac{1}{2}, -\frac{1}{2}, \frac{1}{2})$ & $111111$ & $s_1s_2s_3s_4s_5s_6$ \\
		\hline
		$(0, 0, 0, 1, 1, 0, 0, 0)$ & $011221$ & $s_5s_4s_2s_3s_4s_5s_6$ \\
		\hline
		$(\frac{1}{2}, -\frac{1}{2}, \frac{1}{2}, -\frac{1}{2}, \frac{1}{2}, -\frac{1}{2}, -\frac{1}{2}, \frac{1}{2})$ & $111211$ & $s_4s_1s_2s_3s_4s_5s_6$ \\
		\hline
		$(-\frac{1}{2}, \frac{1}{2}, \frac{1}{2}, -\frac{1}{2}, \frac{1}{2}, -\frac{1}{2}, -\frac{1}{2}, \frac{1}{2})$ & $112211$ & $s_3s_4s_1s_2s_3s_4s_5s_6$ \\
		\hline
		$(\frac{1}{2}, -\frac{1}{2}, -\frac{1}{2}, \frac{1}{2}, \frac{1}{2}, -\frac{1}{2}, -\frac{1}{2}, \frac{1}{2})$ & $111221$ & $s_5s_4s_1s_2s_3s_4s_5s_6$ \\
		\hline
		$(-\frac{1}{2}, \frac{1}{2}, -\frac{1}{2}, \frac{1}{2}, \frac{1}{2}, -\frac{1}{2}, -\frac{1}{2}, \frac{1}{2})$ & $112221$ & $s_3s_5s_4s_1s_2s_3s_4s_5s_6$ \\
		\hline
		$(-\frac{1}{2}, -\frac{1}{2}, \frac{1}{2}, \frac{1}{2}, \frac{1}{2}, -\frac{1}{2}, -\frac{1}{2}, \frac{1}{2})$ & $112321$ & $s_4s_3s_5s_4s_1s_2s_3s_4s_5s_6$ \\
		\hline
		$(\frac{1}{2}, \frac{1}{2}, \frac{1}{2}, \frac{1}{2}, \frac{1}{2}, -\frac{1}{2}, -\frac{1}{2}, \frac{1}{2})$ & $122321$ & $s_2s_4s_3s_5s_4s_1s_2s_3s_4s_5s_6$ \\
		\hline
	\end{longtable}

\begin{longtable}{|c|c|c|}
		\hline
\multicolumn{3}{|c|}{Case $\Phi=E_7$ with the order $\alpha_3,~\alpha_7,~\alpha_4,~\alpha_6,~\alpha_5,~\alpha_2,~\alpha_1$}\\\hline
		$(c_1,~\ldots,~c_8)$ & $(b_1,~\ldots,~b_7)$& Reduced decomposition of $u$\\
		\hline\hline
		$(0, 0, 0, 0, -1, 1, 0, 0)$ & $0000001$ & $s_7$ \\
		\hline
		$(0, 0, 0, -1, 0, 1, 0, 0)$ & $0000011$ & $s_6s_7$ \\
		\hline
		$(0, 0, -1, 0, 0, 1, 0, 0)$ & $0000111$ & $s_5s_6s_7$\\
		\hline
		$(0, -1, 0, 0, 0, 1, 0, 0)$ & $0001111$ & $s_4s_5s_6s_7$\\
		\hline
		$(1, 0, 0, 0, 0, 1, 0, 0)$ & $0101111$ & $s_2s_4s_5s_6s_7$\\
		\hline
		$(-1, 0, 0, 0, 0, 1, 0, 0)$ & $0011111$ & $s_3s_4s_5s_6s_7$ \\
		\hline
		$(0, 1, 0, 0, 0, 1, 0, 0)$ & $0111111$ & $s_2s_3s_4s_5s_6s_7$\\
		\hline
		$(-\frac{1}{2}, -\frac{1}{2}, -\frac{1}{2}, -\frac{1}{2}, -\frac{1}{2}, \frac{1}{2}, -\frac{1}{2}, \frac{1}{2})$ & $1011111$ & $s_1s_3s_4s_5s_6s_7$\\
		\hline
		$(\frac{1}{2}, \frac{1}{2}, -\frac{1}{2}, -\frac{1}{2}, -\frac{1}{2}, \frac{1}{2}, -\frac{1}{2}, \frac{1}{2})$ & $1111111$ & $s_2s_1s_3s_4s_5s_6s_7$\\
		\hline
		$(0, 0, 1, 0, 0, 1, 0, 0)$ & $0112111$ &  $s_4s_2s_3s_4s_5s_6s_7$\\
		\hline
		$(\frac{1}{2}, -\frac{1}{2}, \frac{1}{2}, -\frac{1}{2}, -\frac{1}{2}, \frac{1}{2}, -\frac{1}{2}, \frac{1}{2})$ & $1112111$ & $s_1s_4s_2s_3s_4s_5s_6s_7$\\
		\hline
		$(0, 0, 0, 1, 0, 1, 0, 0)$ & $0112211$ & $s_5s_4s_2s_3s_4s_5s_6s_7$\\
		\hline
		$(\frac{1}{2}, -\frac{1}{2}, -\frac{1}{2}, \frac{1}{2}, -\frac{1}{2}, \frac{1}{2}, -\frac{1}{2}, \frac{1}{2})$ & $1112211$ & $s_5s_1s_4s_2s_3s_4s_5s_6s_7$\\
		\hline
		$(0, 0, 0, 0, 1, 1, 0, 0)$ & $0112221$ & $s_6s_5s_4s_2s_3s_4s_5s_6s_7$\\
		\hline
		$(-\frac{1}{2}, \frac{1}{2}, \frac{1}{2}, -\frac{1}{2}, -\frac{1}{2}, \frac{1}{2}, -\frac{1}{2}, \frac{1}{2})$ & $1122111$ & $s_3s_1s_4s_2s_3s_4s_5s_6s_7$\\
		\hline
		$(\frac{1}{2}, -\frac{1}{2}, -\frac{1}{2}, -\frac{1}{2}, \frac{1}{2}, \frac{1}{2}, -\frac{1}{2}, \frac{1}{2})$ & $1112221$ & $s_6s_5s_1s_4s_2s_3s_4s_5s_6s_7$ \\
		\hline
		$(-\frac{1}{2}, \frac{1}{2}, -\frac{1}{2}, \frac{1}{2}, -\frac{1}{2}, \frac{1}{2}, -\frac{1}{2}, \frac{1}{2})$ & $1122211$ & $s_3s_5s_1s_4s_2s_3s_4s_5s_6s_7$ \\
		\hline
		$(-\frac{1}{2}, -\frac{1}{2}, \frac{1}{2}, \frac{1}{2}, -\frac{1}{2}, \frac{1}{2}, -\frac{1}{2}, \frac{1}{2})$ & $1123211$ & $s_4s_3s_5s_1s_4s_2s_3s_4s_5s_6s_7$ \\
		\hline
		$(-\frac{1}{2}, \frac{1}{2}, -\frac{1}{2}, -\frac{1}{2}, \frac{1}{2}, \frac{1}{2}, -\frac{1}{2}, \frac{1}{2})$ & $1122221$ & $s_6s_3s_5s_1s_4s_2s_3s_4s_5s_6s_7$ \\
		\hline
		$(-\frac{1}{2}, -\frac{1}{2}, \frac{1}{2}, -\frac{1}{2}, \frac{1}{2}, \frac{1}{2}, -\frac{1}{2}, \frac{1}{2})$ & $1123221$ & $s_4s_3s_5s_1s_4s_2s_3s_4s_5s_6s_7$ \\
		\hline
		$(\frac{1}{2}, \frac{1}{2}, \frac{1}{2}, \frac{1}{2}, -\frac{1}{2}, \frac{1}{2}, -\frac{1}{2}, \frac{1}{2})$ & $1223211$ & $s_2s_4s_3s_5s_1s_4s_2s_3s_4s_5s_6s_7$\\
		\hline
		$(-\frac{1}{2}, -\frac{1}{2}, -\frac{1}{2}, \frac{1}{2}, \frac{1}{2}, \frac{1}{2}, -\frac{1}{2}, \frac{1}{2})$ & $1123321$ & $s_5s_4s_3s_5s_1s_4s_2s_3s_4s_5s_6s_7$ \\
		\hline
		$(\frac{1}{2}, \frac{1}{2}, \frac{1}{2}, -\frac{1}{2}, \frac{1}{2}, \frac{1}{2}, -\frac{1}{2}, \frac{1}{2})$ & $1223221$ & $s_6s_2s_4s_3s_5s_1$\\ & &$s_4s_2s_3s_4s_5s_6s_7$ \\
		\hline
		$(\frac{1}{2}, \frac{1}{2}, -\frac{1}{2}, \frac{1}{2}, \frac{1}{2}, \frac{1}{2}, -\frac{1}{2}, \frac{1}{2})$ & $1223321$ & $s_5s_6s_2s_4s_3s_5s_1$\\ & & $s_4s_2s_3s_4s_5s_6s_7$ \\
		\hline
		$(\frac{1}{2}, -\frac{1}{2}, \frac{1}{2}, \frac{1}{2}, \frac{1}{2}, \frac{1}{2}, -\frac{1}{2}, \frac{1}{2})$ & $1224321$ & $s_4s_5s_6s_2s_4s_3s_5s_1$ \\ & & $s_4s_2s_3s_4s_5s_6s_7$ \\
		\hline
		$(-\frac{1}{2}, \frac{1}{2}, \frac{1}{2}, \frac{1}{2}, \frac{1}{2}, \frac{1}{2}, -\frac{1}{2}, \frac{1}{2})$ & $1234321$ & $s_3s_4s_5s_6s_2s_4s_3s_5s_1$ \\ & &$s_4s_2s_3s_4s_5s_6s_7$ \\
		\hline
		$(0, 0, 0, 0, 0, 0, -1, 1)$ & $2234321$ & $s_1s_3s_4s_5s_6s_2s_4s_3s_5s_1$ \\ & & $s_4s_2s_3s_4s_5s_6s_7$ \\
		\hline
	\end{longtable}

\begin{longtable}{|c|c|c|}
		\hline
\multicolumn{3}{|c|}{Case $\Phi=E_8$ with the order $\alpha_4,~\alpha_8,~\alpha_5,~\alpha_7,~\alpha_6,~\alpha_3,~\alpha_1,~\alpha_1$}\\\hline
		$(c_1,~\ldots,~c_8)$ & $(b_1,~\ldots,~b_7)$& Reduced decomposition of $u$\\
		\hline\hline
		$(0, 0, 0, 0, 0, -1, 1, 0)$ & $00000001$ & $s_8$\\
		\hline
		$(0, 0, 0, 0, -1, 0, 1, 0)$ & $00000011$ & $s_7s_8$ \\
		\hline
		$(0, 0, 0, -1, 0, 0, 1, 0)$ & $00000111$ &  $s_6s_7s_8$ \\
		\hline
		$(0, 0, -1, 0, 0, 0, 1, 0)$ & $00001111$ & $s_5s_6s_7s_8$ \\
		\hline
		$(0, -1, 0, 0, 0, 0, 1, 0)$ & $00011111$ & $s_4s_5s_6s_7s_8$\\
		\hline
		$(-1, 0, 0, 0, 0, 0, 1, 0)$ & $00111111$ & $s_3s_4s_5s_6s_7s_8$\\
		\hline
		$(1, 0, 0, 0, 0, 0, 1, 0)$ & $01011111$ & $s_2s_4s_5s_6s_7s_8$\\
		\hline
		$(-\frac{1}{2}, -\frac{1}{2}, -\frac{1}{2}, -\frac{1}{2}, -\frac{1}{2}, -\frac{1}{2}, \frac{1}{2}, \frac{1}{2})$ & $10111111$ & $s_1s_3s_4s_5s_6s_7s_8$\\
		\hline
		$(0, 1, 0, 0, 0, 0, 1, 0)$ & $01111111$ & $s_3s_2s_4s_5s_6s_7s_8$\\
		\hline
		$(0, 0, 1, 0, 0, 0, 1, 0)$ & $01121111$ & $s_4s_3s_2s_4s_5s_6s_7s_8$\\
		\hline
		$(\frac{1}{2}, \frac{1}{2}, -\frac{1}{2}, -\frac{1}{2}, -\frac{1}{2}, -\frac{1}{2}, \frac{1}{2}, \frac{1}{2})$ & $11111111$ & $s_1s_3s_2s_4s_5s_6s_7s_8$\\
		\hline
		$(0, 0, 0, 1, 0, 0, 1, 0)$ & $01122111$ & $s_5s_4s_3s_2s_4s_5s_6s_7s_8$\\
		\hline
		$(\frac{1}{2}, -\frac{1}{2}, \frac{1}{2}, -\frac{1}{2}, -\frac{1}{2}, -\frac{1}{2}, \frac{1}{2}, \frac{1}{2})$ & $11121111$ & $s_1s_4s_3s_2s_4s_5s_6s_7s_8$\\
		\hline
		$(\frac{1}{2}, -\frac{1}{2}, -\frac{1}{2}, \frac{1}{2}, -\frac{1}{2}, -\frac{1}{2}, \frac{1}{2}, \frac{1}{2})$ & $11122111$ & $s_5s_1s_4s_3s_2s_4s_5s_6s_7s_8$\\
		\hline
		$(-\frac{1}{2}, \frac{1}{2}, \frac{1}{2}, -\frac{1}{2}, -\frac{1}{2}, -\frac{1}{2}, \frac{1}{2}, \frac{1}{2})$ & $11221111$ & $s_3s_1s_4s_3s_2s_4s_5s_6s_7s_8$\\
		\hline
		$(0, 0, 0, 0, 1, 0, 1, 0)$ & $01122211$ & $s_6s_5s_1s_4s_3s_2s_4s_5s_6s_7s_8$\\
		\hline
		$(-\frac{1}{2}, \frac{1}{2}, -\frac{1}{2}, \frac{1}{2}, -\frac{1}{2}, -\frac{1}{2}, \frac{1}{2}, \frac{1}{2})$ & $11222111$ & $s_5s_3s_1s_4s_3s_2s_4s_5s_6s_7s_8$\\
		\hline
		$(0, 0, 0, 0, 0, 1, 1, 0)$ & $01122221$ & $s_7s_6s_5s_1s_4s_3s_2s_4s_5s_6s_7s_8$\\
		\hline
		$(\frac{1}{2}, -\frac{1}{2}, -\frac{1}{2}, -\frac{1}{2}, \frac{1}{2}, -\frac{1}{2}, \frac{1}{2}, \frac{1}{2})$ & $11122211$ & $s_1s_6s_5s_1s_4s_3s_2$ \\ & & $s_4s_5s_6s_7s_8$\\
		\hline
		$(\frac{1}{2}, -\frac{1}{2}, -\frac{1}{2}, -\frac{1}{2}, -\frac{1}{2}, \frac{1}{2}, \frac{1}{2}, \frac{1}{2})$ & $11122221$ & $s_7s_1s_6s_5s_1s_4s_3$ \\ & & $s_2s_4s_5s_6s_7s_8$\\
		\hline
		$(-\frac{1}{2}, \frac{1}{2}, -\frac{1}{2}, -\frac{1}{2}, \frac{1}{2}, -\frac{1}{2}, \frac{1}{2}, \frac{1}{2})$ & $11222211$ & $s_3s_1s_6s_5s_1s_4s_3$ \\ & & $s_2s_4s_5s_6s_7s_8$\\
		\hline
		$(-\frac{1}{2}, -\frac{1}{2}, \frac{1}{2}, \frac{1}{2}, -\frac{1}{2}, -\frac{1}{2}, \frac{1}{2}, \frac{1}{2})$ & $11232111$ & $s_4s_5s_3s_1s_4s_3$ \\ & & $s_2s_4s_5s_6s_7s_8$\\
		\hline
		$(-\frac{1}{2}, -\frac{1}{2}, \frac{1}{2}, -\frac{1}{2}, \frac{1}{2}, -\frac{1}{2}, \frac{1}{2}, \frac{1}{2})$ & $11232211$ & $s_6s_4s_5s_3s_1s_4$ \\ & & $s_3s_2s_4s_5s_6s_7s_8$\\
		\hline
		$(-\frac{1}{2}, \frac{1}{2}, -\frac{1}{2}, -\frac{1}{2}, -\frac{1}{2}, \frac{1}{2}, \frac{1}{2}, \frac{1}{2})$ & $11222221$ & $s_7s_3s_1s_6s_5s_1s_4s_3$ \\ & & $s_2s_4s_5s_6s_7s_8$\\
		\hline
		$(\frac{1}{2}, \frac{1}{2}, \frac{1}{2}, \frac{1}{2}, -\frac{1}{2}, -\frac{1}{2}, \frac{1}{2}, \frac{1}{2})$ & $12232111$ & $s_2s_4s_5s_3s_1s_4s_3s_2$ \\ & & $s_4s_5s_6s_7s_8$\\
		\hline
		$(-\frac{1}{2}, -\frac{1}{2}, -\frac{1}{2}, \frac{1}{2}, \frac{1}{2}, -\frac{1}{2}, \frac{1}{2}, \frac{1}{2})$ & $11233211$ & $s_5s_6s_4s_5s_3s_1s_4s_3$ \\ & & $s_2s_4s_5s_6s_7s_8$\\
		\hline
		$(\frac{1}{2}, \frac{1}{2}, \frac{1}{2}, -\frac{1}{2}, \frac{1}{2}, -\frac{1}{2}, \frac{1}{2}, \frac{1}{2})$ & $12232211$ & $s_6s_2s_4s_5s_3s_1s_4s_3$ \\ & & $s_2s_4s_5s_6s_7s_8$\\
		\hline
		$(-\frac{1}{2}, -\frac{1}{2}, \frac{1}{2}, -\frac{1}{2}, -\frac{1}{2}, \frac{1}{2}, \frac{1}{2}, \frac{1}{2})$ & $11232221$ & $s_7s_6s_4s_5s_3s_1s_4s_3$ \\ & & $s_2s_4s_5s_6s_7s_8$\\
		\hline
		$(\frac{1}{2}, \frac{1}{2}, -\frac{1}{2}, \frac{1}{2}, \frac{1}{2}, -\frac{1}{2}, \frac{1}{2}, \frac{1}{2})$ & $12233211$ & $s_5s_6s_2s_4s_5s_3s_1s_4$ \\ & & $s_3s_2s_4s_5s_6s_7s_8$\\
		\hline
		$(\frac{1}{2}, \frac{1}{2}, \frac{1}{2}, -\frac{1}{2}, -\frac{1}{2}, \frac{1}{2}, \frac{1}{2}, \frac{1}{2})$ & $12232221$ & $s_2s_7s_6s_4s_5s_3s_1s_4$ \\ & & $s_3s_2s_4s_5s_6s_7s_8$\\
		\hline
		$(-\frac{1}{2}, -\frac{1}{2}, -\frac{1}{2}, \frac{1}{2}, -\frac{1}{2}, \frac{1}{2}, \frac{1}{2}, \frac{1}{2})$ & $11233221$ & $s_7s_5s_6s_4s_5s_3s_1s_4$ \\ & & $s_3s_2s_4s_5s_6s_7s_8$\\
		\hline
		$(-\frac{1}{2}, -\frac{1}{2}, -\frac{1}{2}, -\frac{1}{2}, \frac{1}{2}, \frac{1}{2}, \frac{1}{2}, \frac{1}{2})$ & $11233321$ & $s_6s_7s_5s_6s_4s_5s_3s_1$ \\ & & $s_4s_3s_2s_4s_5s_6s_7s_8$\\
		\hline
		$(\frac{1}{2}, \frac{1}{2}, -\frac{1}{2}, \frac{1}{2}, -\frac{1}{2}, \frac{1}{2}, \frac{1}{2}, \frac{1}{2})$ & $12233221$ & $s_5s_2s_7s_6s_4s_5s_3s_1$ \\ & & $s_4s_3s_2s_4s_5s_6s_7s_8$\\
		\hline
		$(\frac{1}{2}, -\frac{1}{2}, \frac{1}{2}, \frac{1}{2}, \frac{1}{2}, -\frac{1}{2}, \frac{1}{2}, \frac{1}{2})$ & $12243211$ & $s_4s_5s_6s_2s_4s_5s_3s_1$ \\ & & $s_4s_3s_2s_4s_5s_6s_7s_8$\\
		\hline
		$(-\frac{1}{2}, \frac{1}{2}, \frac{1}{2}, \frac{1}{2}, \frac{1}{2}, -\frac{1}{2}, \frac{1}{2}, \frac{1}{2})$ & $12343211$ & $s_3s_4s_5s_6s_2s_4s_5s_3s_1$ \\ & & $s_4s_3s_2s_4s_5s_6s_7s_8$\\
		\hline
		$(\frac{1}{2}, -\frac{1}{2}, \frac{1}{2}, \frac{1}{2}, -\frac{1}{2}, \frac{1}{2}, \frac{1}{2}, \frac{1}{2})$ & $12243221$ & $s_7s_4s_5s_6s_2s_4s_5s_3s_1$ \\ & & $s_4s_3s_2s_4s_5s_6s_7s_8$ \\
		\hline
		$(\frac{1}{2}, \frac{1}{2}, -\frac{1}{2}, -\frac{1}{2}, \frac{1}{2}, \frac{1}{2}, \frac{1}{2}, \frac{1}{2})$ & $12233321$ & $s_6s_5s_2s_7s_6s_4s_5s_3s_1$ \\ & & $s_4s_3s_2s_4s_5s_6s_7s_8$\\
		\hline
		$(\frac{1}{2}, -\frac{1}{2}, \frac{1}{2}, -\frac{1}{2}, \frac{1}{2}, \frac{1}{2}, \frac{1}{2}, \frac{1}{2})$ & $12243321$ & $s_4s_6s_5s_2s_7s_6s_4s_5s_3$ \\ & & $s_1s_4s_3s_2s_4s_5s_6s_7s_8$\\
		\hline
		$(0, 0, 0, 0, 0, -1, 0, 1)$ & $22343211$ & $s_1s_3s_4s_5s_6s_2s_4s_5$ \\ & & $s_3s_1s_4s_3s_2s_4s_5s_6s_7s_8$\\
		\hline
		$(-\frac{1}{2}, \frac{1}{2}, \frac{1}{2}, \frac{1}{2}, -\frac{1}{2}, \frac{1}{2}, \frac{1}{2}, \frac{1}{2})$ & $12343221$ & $s_7s_3s_4s_5s_6s_2s_4s_5$ \\ & & $s_3s_1s_4s_3s_2s_4s_5s_6s_7s_8$\\
		\hline
		$(0, 0, 0, 0, -1, 0, 0, 1)$ & $22343221$ & $s_1s_7s_3s_4s_5s_6s_2s_4s_5$ \\ & & $s_3s_1s_4s_3s_2s_4s_5s_6s_7s_8$\\
		\hline
		$(\frac{1}{2}, -\frac{1}{2}, -\frac{1}{2}, \frac{1}{2}, \frac{1}{2}, \frac{1}{2}, \frac{1}{2}, \frac{1}{2})$ & $12244321$ & $s_5s_4s_6s_5s_2s_7s_6s_4s_5$ \\ & & $s_3s_1s_4s_3s_2s_4s_5s_6s_7s_8$\\
		\hline
		$(-\frac{1}{2}, \frac{1}{2}, \frac{1}{2}, -\frac{1}{2}, \frac{1}{2}, \frac{1}{2}, \frac{1}{2}, \frac{1}{2})$ & $12343321$ & $s_6s_7s_3s_4s_5s_6s_2s_4s_5$ \\ & & $s_3s_1s_4s_3s_2s_4s_5s_6s_7s_8$\\
		\hline
		$(0, 0, 0, -1, 0, 0, 0, 1)$ & $22343321$ & $s_1s_6s_7s_3s_4s_5s_6s_2s_4s_5$ \\ & & $s_3s_1s_4s_3s_2s_4s_5s_6s_7s_8$\\
		\hline
		$(-\frac{1}{2}, \frac{1}{2}, -\frac{1}{2}, \frac{1}{2}, \frac{1}{2}, \frac{1}{2}, \frac{1}{2}, \frac{1}{2})$ & $12344321$ & $s_5s_6s_7s_3s_4s_5s_6s_2s_4s_5$ \\ & & $s_3s_1s_4s_3s_2s_4s_5s_6s_7s_8$\\
		\hline
		$(0, 0, -1, 0, 0, 0, 0, 1)$ & $22344321$ & $s_1s_5s_6s_7s_3s_4s_5s_6s_2s_4$ \\ & & $s_5s_3s_1s_4s_3s_2s_4s_5s_6s_7s_8$\\
		\hline
		$(-\frac{1}{2}, -\frac{1}{2}, \frac{1}{2}, \frac{1}{2}, \frac{1}{2}, \frac{1}{2}, \frac{1}{2}, \frac{1}{2})$ & $12354321$ & $s_4s_5s_6s_7s_3s_4s_5s_6s_2s_4$ \\ & & $s_5s_3s_1s_4s_3s_2s_4s_5s_6s_7s_8$\\
		\hline
		$(\frac{1}{2}, \frac{1}{2}, \frac{1}{2}, \frac{1}{2}, \frac{1}{2}, \frac{1}{2}, \frac{1}{2}, \frac{1}{2})$ & $13354321$ & $s_2s_4s_5s_6s_7s_3s_4s_5s_6s_2s_4$ \\ & & $s_5s_3s_1s_4s_3s_2s_4s_5s_6s_7s_8$\\
		\hline
		$(0, -1, 0, 0, 0, 0, 0, 1)$ & $22354321$ & $s_1s_4s_5s_6s_7s_3s_4s_5s_6s_2s_4$ \\ & & $s_5s_3s_1s_4s_3s_2s_4s_5s_6s_7s_8$\\
		\hline
		$(1, 0, 0, 0, 0, 0, 0, 1)$ & $23354321$ & $s_2s_1s_4s_5s_6s_7s_3s_4s_5s_6s_2s_4$ \\ & & $s_5s_3s_1s_4s_3s_2s_4s_5s_6s_7s_8$\\
		\hline
		$(-1, 0, 0, 0, 0, 0, 0, 1)$ & $22454321$ & $s_3s_1s_4s_5s_6s_7s_3s_4s_5s_6s_2s_4$ \\ & & $s_5s_3s_1s_4s_3s_2s_4s_5s_6s_7s_8$\\
		\hline
		$(0, 1, 0, 0, 0, 0, 0, 1)$ & $23454321$ & $s_2s_3s_1s_4s_5s_6s_7s_3s_4s_5s_6s_2s_4$ \\ & & $s_5s_3s_1s_4s_3s_2s_4s_5s_6s_7s_8$\\
		\hline
		$(0, 0, 1, 0, 0, 0, 0, 1)$ & $23464321$ & $s_4s_2s_3s_1s_4s_5s_6$ \\ & & $s_7s_3s_4s_5s_6s_2s_4s_5$ \\ & & $s_3s_1s_4s_3s_2s_4s_5s_6s_7s_8$\\
		\hline
		$(0, 0, 0, 1, 0, 0, 0, 1)$ & $23465321$ & $s_5s_4s_2s_3s_1s_4s_5s_6$ \\ & & $s_7s_3s_4s_5s_6s_2s_4$ \\ & & $s_5s_3s_1s_4s_3s_2s_4s_5s_6s_7s_8$\\
		\hline
		$(0, 0, 0, 0, 1, 0, 0, 1)$ & $23465421$ & $s_6s_5s_4s_2s_3s_1s_4s_5$ \\ & & $s_6s_7s_3s_4s_5s_6$ \\ & & $s_2s_4s_5s_3s_1s_4s_3s_2s_4s_5s_6s_7s_8$\\
		\hline
		$(0, 0, 0, 0, 0, 1, 0, 1)$ & $23465431$ & $s_7s_6s_5s_4s_2s_3s_1s_4$ \\ & & $s_5s_6s_7s_3s_4s_5s_6$ \\ & & $s_2s_4s_5s_3s_1s_4s_3s_2s_4s_5s_6s_7s_8$\\
		\hline	
	\end{longtable}

All these tables were generated using computer algebra system SAGE \cite{SAGE}; the listing of the code can be found in the Appendix.

One can immediately check that $s_{\beta}$ and $u$ satisfy the conditions of Proposition~\ref{prop:crucial}. Hence, according to this proposition, $\beta_i$ does not divide $d_{w_i}$ in $\Cp[\htt]$ for $i=1,~2$. On the other hand, formula (\ref{formula:dyer}) implies that, given $w\in W$, there exists $g\in\Cp[\htt]$ such that
\begin{equation*}
d_w=c_w\cdot\prod_{\alpha\in\Phi^+}\alpha=g\cdot\prod_{\alpha\in\Phi^+,~s_{\alpha}\nleq w}\alpha,
\end{equation*}
hence if $s_{\alpha}\neq w$ then $\alpha$ divides $d_w$. But if, for example, $s_{\beta_1}\nleq w_2$ then $\beta_1$ divides $d_{w_2}$. At the same time, $\beta_1$ does not divide $d_{w_1}$, thus, $d_{w_1}\neq d_{w_2}$. The proof is complete.  }

\medskip\textsc{Mikhail V. Ignatyev: Samara National Research University, Ak. Pavlova 1, 443011,\\\indent Samara, Russia}

\emph{E-mail address}: \texttt{mihail.ignatev@gmail.com}

\medskip\textsc{Aleksandr A. Shevchenko: Samara National Research University, Ak. Pavlova 1,\\\indent 443011, Samara, Russia}

\emph{E-mail address}: \texttt{shevchenko.alexander.1618@gmail.com }

\newpage
\begin{center}
\textbf{Appendix}
\end{center}

Below we present the listing of the code generating tables from the proof of Proposition~\ref{prop:good_pair} using computer algebra system SAGE.

\begin{verbatim}
rank=8 # the rank of the root system
column_number=8 # the number of the first column
W=WeylGroup(['E',rank],prefix='s', implementation='permutation')
ref=W.reflections()
s=W.simple_reflections()
R=RootSystem(['E',rank]).ambient_space();
simple_roots=R.simple_roots()
phi_plus=W.positive_roots()
C1=[]
C1el=[]
for i in range(0,len(phi_plus)):
  if phi_plus[i][column_number-1]!=0:
    C1.append(phi_plus[i])
    C1el.append(ref[i+1])
U=[s[column_number]]
Ulistver=[[column_number]]
for i in range(1,len(C1)):
  u=copy(Ulistver[i-1])
  index=-1
  difference=C1[i]-C1[i-1]
  difference_abs=[abs(ele) for ele in difference]
  if sum(difference_abs)!=1:
    b=0
    for j in range(i-2,1,-1):
      difference1=C1[i]-C1[j]
      difference_abs1=[abs(ele) for ele in difference1]
      if sum(difference_abs1)==1 and b==0:
        b=1
        difference=C1[i]-C1[j]
        difference_abs=[abs(ele) for ele in difference]
        u=copy(Ulistver[j])
  if sum(difference)==1:
    for j in range(0,len(difference)):
      if difference[j]==1:
        index=j
  u1=[]
  u2=s[column_number]*s[column_number]
  if index!=-1:
    u1=[index+1]
    u2=s[index+1]
    for j in range(0,len(u)):
      u1.append(u[j])
      u2=u2*s[u[j]]
  U.append(u2)
  Ulistver.append(u1)
list_of_indexes=[]
for i in range(1,rank+1):
  if i!=column_number:
    list_of_indexes.append(i)
for u in U:
  b=1
  for i in list_of_indexes:
    u1=u*s[i]
    if u1.length()<u.length():
      b=0
    if b==0:
      print('false')
      print(u)
V=[s[column_number]*s[column_number]]
Vlistver=[[]]
for i in range(1,len(Ulistver)):
  u=copy(Ulistver[i])
  v1=[]
  v2=s[column_number]*s[column_number]
  for j in range(len(u)-2,-1,-1):
    v1.append(u[j])
    v2=v2*s[u[j]]
  V.append(v2)
  Vlistver.append(v1)
C1roots=[R.simple_root(column_number)]
for i in range(1,len(Vlistver)):
  root=R.simple_root(column_number)
  for j in range(0,len(Vlistver[i])):
    root=root+R.simple_root(Vlistver[i][j])
  C1roots.append(root)
for i in range(0,len(U)):
  u=U[i]
  v=V[i]
  r=C1el[i]
  if (u*v!=r) or (u.length()+v.length()!=r.length()):
    print(false)
for i in range(0,len(U)):
  print(C1roots[i],C1[i],U[i])
\end{verbatim}

\end{document}